\documentclass[12pt,fleqn]{amsart}

\usepackage[all]{xy}
\usepackage{tikz-cd}
\usepackage{amscd}
\usepackage{color}
\usepackage{enumerate}
\usepackage{amssymb}
\newtheorem{theorem}{Theorem}[section]
\newtheorem{lemma}[theorem]{Lemma}
\newtheorem{proposition}[theorem]{Proposition}
\newtheorem{problem}[theorem]{Problem}

\theoremstyle{definition}
\newtheorem{definition}[theorem]{Definition}

\newtheorem{remark}[theorem]{Remark}

\newtheorem*{ind}{Key Property KP($\xi$)}
\theoremstyle{remark}

\numberwithin{equation}{section}

\newcommand{\nnorm}[1]{{\left\vert\kern-0.25ex\left\vert\kern-0.25ex\left\vert #1
		\right\vert\kern-0.25ex\right\vert\kern-0.25ex\right\vert}}

\newcommand{\cA}{\mathcal A}
\newcommand{\cB}{\mathcal B}
\newcommand{\cC}{\mathcal C}

\newcommand{\cP}{\mathcal P}

\newcommand{\cZ}{\mathcal Z}

\newcommand{\R}{\mathbb R}
\newcommand{\cR}{\mathcal R}

\newcommand{\fA}{\mathfrak A}
\newcommand{\fB}{\mathfrak B}

\newcommand{\cF}{\mathcal F}

\newcommand{\ult}{\protect{\rm ult}}

\newcommand{\sub}{\subseteq}
\newcommand{\eps}{\varepsilon}
\newcommand{\er}{\mathbb R}
\newcommand{\sm}{\setminus}
\newcommand{\con}{\mathfrak c}

\newcommand{\vf}{\varphi}
\newcommand{\fin}{\protect{\rm fin }}

\newcommand{\ml}{M^{\ell_1}}
\newcommand{\jl}{\protect{\rm JL}}
\newcommand{\bil}[2]{{\left\langle\kern0ex #1,#2
		\kern0ex\right\rangle}}

\begin{document}


\baselineskip=17pt

 \title[The  complemented subspace problem]
 {The  complemented subspace problem for    $C(K)$-spaces:\\ A counterexample}

\author[G.\ Plebanek]{Grzegorz Plebanek}
\address{Instytut Matematyczny\\ Uniwersytet Wroc\l awski\\ Pl.\ Grunwaldzki 2/4\\
50-384 Wroc\-\l aw\\ Poland} \email{grzegorz.plebanek@math.uni.wroc.pl}

\author[A. Salguero Alarc\'on]{Alberto Salguero Alarc\'on}
\address{Instituto de Matem\'aticas\\ Universidad de Extremadura\\
Avenida de Elvas\\ 06071-Badajoz\\ Spain} \email{salgueroalarcon@unex.es}

\thanks{The first author was partially supported by the grant
2018/29/B/ST1/00223 from National Science Centre, Poland.
The second author  has been supported  by project PID2019-103961GB-C21 funded by MCIN.
project IB20038 funded by Junta de Extremadura and by the grant FPU18/00990 by the
Spanish Ministry of Universities.
}


\begin{abstract}
We construct a compact space $L$ and a $1$-complemented subspace of the Banach space $C(L)$
which is not isomorphic to a space of continuous functions.
\end{abstract}

\date{}
\maketitle

\section{Introduction}
Let us say that a Banach space $X$ is a \emph{$C$-space} if
$X$ is isomorphic to a space of the form $C(K)$, the space of all real-valued continuous functions on
a Hausdorff compact space $K$.
The following complemented subspace problem for $C$-spaces
has been around for years, see Lindenstrauss and Wulbert \cite[p.\ 348]{LW69}:
\begin{problem}\label{csp}
	Is the class of  $C$-spaces closed under taking complemented subspaces?
\end{problem}

Recall that a hyperplane in a $C$-space is again a $C$-space.
Indeed, if $H$ is a hyperplane in the space $C(K)$, then, as all the hyperplanes of a given Banach space are
isomorphic, we can assume that $H=\{g\in C(K)\colon g(x_0)=g(x_1)\}$ for some distinct $x_0,x_1\in K$.
But then $H$ is isomorphic to $C(\widetilde{K})$, where
$\widetilde{K}$ is obtained from $K$ by identifying the points $x_0$ and $x_1$.
Nonetheless,  a hyperplane in a space of the form $C(K)$ need not be isomorphic
to the space $C(K)$ itself. In \cite{Ko04}, Koszmider gives an involved construction of such a peculiar compact space.

 Rosenthal's survey \cite[Section 5]{Ro03} offers nearly everything that is known on the nature
 of complemented subspaces of $C$-spaces.
 The following are examples of $C$-spaces for which a positive answer to Problem \ref{csp} has been obtained:

\begin{enumerate}
	\item A subspace of $c_0(\Gamma)$ is complemented if and only if it is isomorphic to some $c_0(\Gamma')$ with $|\Gamma'|\leq |\Gamma|$, see Granero \cite{Gra98}.
	\item Every complemented subspace of $C(\omega^\omega)$ is isomorphic either to $c_0$ or to $C(\omega^\omega)$, see Benyamini \cite{Ben78}.
	\item Every complemented subspace of $C[0,1]$ with a non-separable dual is isomorphic to $C[0,1]$;
this was proved  by Rosenthal \cite{Ro72}, see also \cite[section 4C]{Ro03}.
	\item Every complemented subspace of $\ell_\infty$ is isomorphic to $\ell_\infty$,
by a well-known result of Lindenstrauss \cite{Lin67}.
	\item Every complemented subspace of the space $\ell_\infty^c(\Gamma)$ of bounded functions $\Gamma \to \R$ with countable support (endowed with the supremum norm) is isomorphic  either to $\ell_\infty$ or to $\ell_\infty^c(\Gamma')$ for $|\Gamma'|\leq |\Gamma|$; this is a result of Johnson, Kania and Schechtman \cite{JKS16}.
\end{enumerate}

Here $c_0(\Gamma)$ denotes the Banach space of all functions $x\colon\Gamma\to \er$ vanishing at $\infty$ (meaning that
the set $\{\gamma\in\Gamma\colon |x(\gamma)|\ge\eps\}$ is finite for every $\eps>0$).
The space $c_0(\Gamma)$ may be treated as a hyperplane in the space of continuous functions on the  one-point
compactification of the discrete space $\Gamma$ and therefore it is a $C$-space itself.
It is worth recalling in connection with (1) that Marciszewski \cite{Ma03} gave a characterization of those
compact spaces $K$ for which $C(K)$ is isomorphic to $c_0(\Gamma)$ for some $\Gamma$.
On the other hand, let us note that even the nature of  complemented subspaces of separable $C$-spaces is not fully understood.
It is proved in \cite{Ben78} that if $X$ is isomorphic to a complemented subspace of $C[0,1]$
and $X^\ast$ is separable then there is  a countable ordinal number $\alpha$ such that $X$ is a quotient
of  $C(\alpha)$, the space of continuous functions on
the ordinal interval $\{\beta\colon\beta\le\alpha\}$.

One can  reformulate Problem \ref{csp}, in the spirit of Pe\l czy\'nski's dissertation \cite{Pe68}, to the following.

\begin{problem} \label{quotient} Given a continuous surjection $\theta\colon L \to K$ between compact spaces
such that $\theta^\circ[C(K)]$ is a complemented subspace of $C(L)$,
is the Banach space $C(L)/\theta^\circ[C(K)]$ a $C$-space?
\end{problem}

\noindent Here $\theta^\circ\colon C(K)\to C(L)$ stands for the isometric embedding given by $\theta^\circ(g)=g\circ\theta$ for $g\in C(K)$.
 The aim of this paper is to present a negative answer to Problem \ref{csp} by providing
 the following counterexample to Problem \ref{quotient}.

\begin{theorem}\label{main}
There are two separable scattered compacta $K$ and $L$ with the third derivative empty,
a continuous surjection $\theta\colon L \to K$ and a closed subspace $X$ of $C(L)$ such that
	
\begin{enumerate}[(i)]
\item $C(L) \simeq \theta^\circ[C(K)] \oplus X$;
\item the spaces $ \theta^\circ[C(K)]$ and $X$ are both $1$-complemented  in $C(L)$;
\item the space $X$ is not a $C$-space.
\end{enumerate}
\end{theorem}

In connection with  Theorem \ref{main}(ii), it is worth noting  that if a compact space $L$ is
metrizable and $X$ is $1$-complemented subspace of $C(L)$ then $X$ is a $C$-space, see
Rosenthal \cite[Section 5]{Ro03}.
The compact spaces $K$ and $L$ mentioned in Theorem \ref{main}
are defined by almost disjoint families of subsets of some infinite countable set, so, in other words,
the Banach spaces $C(K)$ and $C(L)$ are isometric to the so called Johnson-Lindenstrauss
spaces, originated in \cite[Example 2]{JL74}.
Those spaces
have  found various applications in Banach space theory, see for instance
\cite{AMP20, CCMPS20,Ko05,KL21, PS21}.

Let us briefly sketch the basic idea that is behind our construction leading to Theorem \ref{main}.
It is easy to check that if a Banach space $X$ is isomorphic to a space of the form $C(K)$, then $B_{X^\ast}$,
the dual closed unit ball equipped with the weak$^\ast$ topology, contains a topological copy of $K$.
We use this observation and employ the fundamental technique of
\cite{PS21}:
in Section \ref{sec:final} we carry out an inductive construction of length $\con$
that builds the space $X$ in question and, at the same time, eliminates all the possible candidates
inside $B_{X^\ast}$ that might arise from an isomorphism between $X$ and some $C$-space.
This inductive process requires no additional set-theoretic axioms thanks to a
counting argument based on a result of Haydon \cite{Ha81}, see Section \ref{haydon}.
In the final section, we apply Theorem \ref{main} to improve our result from \cite{PS21}
on twisted sums of $c_0$ with $C$-spaces.
\medskip

{\bf Acknowledgements.}
We are grateful to Jes\'us Castillo for inspiring this research.
We thank Antonio Avil\'es, Piotr Koszmider and Gonzalo Mart\'{\i}nez-Cervantes
for fruitful conversations.

We are  indebted to the anonymous referees for their very extensive reports and numerous  suggestions
that enabled us to improve
the main result and clarify  several aspects of the presentation.

\section{Preliminaries}\label{p}

\subsection{$C$-spaces}
The letters $K,L$ always stand for  compact Hausdorff spaces.
As we have already mentioned, a Banach space $X$ is a \emph{ $C$-space} if
it is isomorphic to a Banach space of continuous real-valued functions $C(K)$ for some $K$.

 Consider an arbitrary  Banach space $X$. Given any weak$^\ast$ compact subset $K$ of
$B_{X^\ast}$, we say that $K$ is \emph{norming}
if there is $c>0$ such that $\sup_{x^\ast\in K}|x^\ast(x)|\ge c\|x\|$ for every $x\in X$.
Note that in such a case the evaluation map $e(x)(x^\ast)=\bil{x^\ast}{x}$
is an isomorphic embedding 	$X \longrightarrow C(K)$
with  $\|e^{-1}\| \le 1/c$.
	
\begin{definition}
We say that a weak$^\ast$ closed subset $K$ of $B_{X^\ast}$ is \emph{free}
if the above evaluation map $e$ is onto,
that is, for every $f\in C(K)$ there is an $x\in X$ such that $f(x^\ast) = \bil{x^\ast}{x}$ for every $x^\ast\in K$.
\end{definition}

\noindent We recall the following simple observation from \cite[Lemma 2.2]{PS21}.

\begin{lemma}\label{p:1}
For a Banach space $X$, the following are equivalent:
\begin{enumerate}[(i)]
\item $X$ is  a $C$-space;
\item there is a weak$^\ast$ closed subset  $K$ of $B_{X^\ast}$ which is norming and free.
\end{enumerate}
\end{lemma}

\subsection{Johnson-Lindenstrauss spaces}\label{jls}
We denote $\omega=\{0,1,2,...\}$ the set of non-negative integers, and $\fin(\omega)$ is the family of all finite subsets of $\omega$. Given a family $\cF\sub{\mathcal P}(\omega)$, we write $[\cF]$ for the smallest subalgebra of subsets of
$\omega$ containing $\cF$.

\par A family $\cA$ of  infinite subsets of $\omega$ (or any other countable infinite set)
is called \emph{almost disjoint} if $A\cap B\in\fin(\omega)$  for any distinct $A, B\in \cA$.
Given an almost disjoint family $\cA\sub {\mathcal P}(\omega)$, the \emph{Johnson-Lindenstrauss space}
$\jl(\cA)$ is defined as the closure of the subspace of $\ell_\infty$ spanned by the characteristic functions of finite sets, together with the characteristic functions $\chi_A$ for $A\in\cA$ and the constant $1$ function, see \cite[Example 2]{JL74}.  It will be convenient to specify this construction as follows.
Let $\fA = [\cA \cup \fin(\omega)]$ and consider
the linear subspace $s(\fA)$ of $\ell_\infty$ consisting of all simple $\fA$-measurable functions, that is,
 functions of the form $f = \sum_{i\leq n} r_i \cdot \chi_{B_i}$ where $n\in \omega$, $r_i\in \R$ and $B_i \in \fA$. Then it is easy to check the following:

\begin{lemma}\label{jl:1}
	The space  $\jl(\cA)$ is the closure of $s(\fA)$ inside $\ell_\infty$.
\end{lemma}

Given an algebra $\fA$ of subsets of $\omega$, we write $M(\fA)$ for the Banach space of all signed {\em finitely}
additive functions on  $\fA$ of bounded variation.
Recall that the norm of a measure $\mu \in M(\fA)$ is given by
$\|\mu\|=|\mu|(\omega)$, where the variation $|\mu|$ is defined for $A\in\fB$ as
\[|\mu|(A)=\sup_{B\in \fA, \, B\sub A} \big(|\mu(B)| + |\mu(A\setminus B)|\big).\]
In the particular case when $\fA=[\cA\cup\fin(\omega)]$ for a certain almost disjoint family $\cA \subseteq \cP(\omega)$, we can recognize
the dual space $\jl(\cA)^\ast$ as $M(\fA)$, see e.g. \cite[III.2]{DS58}.  Indeed,
every $\mu\in M(\fA)$  defines a bounded linear functional on
$s(\fA)$ by integration of simple functions  and such a functional extends uniquely
to $\jl(\cA)$ by Lemma \ref{jl:1}.
In the sequel, we write $\bil{\mu}{g}$ rather than  $\int  g\;{\rm d}\mu$, and we denote by $M_1(\fA)$ the unit ball in $M(\fA)$ endowed with the weak* topology.

\subsection{Stone spaces}
\label{sec:boole}
We recall here how to change the perspective and pass from
$\jl$-spaces to spaces of continuous functions. Indeed, a quick reference to \cite[Theorem 4.2.5]{AK06} shows
that $\jl$-spaces are $C$-spaces.
To provide a handy construction of the underlying compact spaces, it suffices to apply the classical Stone duality
to our algebras. The reader is referred to
Chapter 3 of Koppelberg's handbook \cite{kopp} for the basic facts concerning the topological version
of the Stone duality.

Given a Boolean algebra $\fB$,  $\ult(\fB)$ denotes its Stone space; that is, the set of all ultrafilters on $\fB$
endowed with the topology having as a base  the collection of sets of the form
$\widehat{B} = \{p\in \ult(\fB)\colon B\in p\}$, for  $B\in \fB$. Then $\ult(\fB)$ is a compact Hausdorff space and the assignment $B\to \widehat{B}$ defines a Boolean
isomorphism between $\fB$ and the algebra of clopen subsets of $\ult(\fB)$.
Conversely, every compact Hausdorff space possessing a base of clopen sets (i.e. a $0$-dimensional space)
can be seen as the Stone space of its algebra of clopen sets. This duality yields a correspondence between morphisms of
 Boolean algebras and their  Stone spaces.
 We will only need the following simple fact, which is a particular case of \cite[Theorem 8.2]{kopp}.

\begin{lemma}\label{stone1}
	If $\fA$ is a subalgebra of a Boolean algebra $\fB$, then $\theta\colon\ult(\fB)\to\ult(\fA)$
	defined by $\theta(p)=p\cap \fA$ for $p\in \ult(\fB)$ is a continuous surjection.
\end{lemma}

The following well-known observation provides a realization of $\jl(\cA)$-spaces as $C$-spaces, see e.g.\
Yost \cite[Example 2]{Yo97}.

\begin{lemma} Let $\cA$ be an almost disjoint family of subsets of $\omega$, and $\fA= [\cA \cup \fin(\omega)]$. Then $\jl(\cA)$ is isometrically isomorphic to $C(\ult(\fA))$.
\end{lemma}
\begin{proof} We simply need to observe that $T_0\chi_A=\chi_{\widehat{A}}$ for $A\in\fA$ defines a linear iso\-me\-tric embedding $T_0$ of $s(\fA)$ into  $C(\ult(\fA))$. Then, $T_0$ extends uniquely to an isometry $T\colon \jl(\cA)\to C(\ult(\fA))$,
which is surjective by the Stone-Weierstrass theorem.
\end{proof}

\subsection{Alexandrov-Urysohn compacta} \label{sec:au}
Bearing in mind that every $\jl$-space  is isometrically isomorphic to some $C(K)$,
it is worth noting some basic topological properties of  the underlying compact space $K$.
Recall that a compact space $K$ is \emph{scattered} if every subspace of $K$ has an isolated point.

Consider an infinite almost disjoint family $\cB$ of infinite subsets of some countable infinite
set $S$ and the algebra $\fB=[\cB\cup\fin(S)]$. Then the space $K=\ult(\fB)$ can be seen
as $S\cup \{p_B\colon B\in \cB\} \cup\{p_\infty\}$, where
\begin{enumerate}
	\item every $s\in S$ is identified with the  \emph{principal} ultrafilter
	$p_s = \{C\in \fB\colon s\in C\}$;
	\item given $B\in \cB$,  $p_B$ is the only ultrafilter containing $B$ but no finite set;
	\item $p_\infty$ is the only ultrafilter on $\fB$ which contains no finite set and no set from $\cB$.
\end{enumerate}

Compact spaces introduced in the above manner  were used nearly a hundred years ago by Alexandroff and Urysohn
and became a popular subject in set-theoretic topology, see Hru\v{s}\'ak's survey \cite{Hr14}.
Locally compact spaces of the form $\ult(\fB)\sm\{p_\infty\}$ are often called  Mr\'owka
or Mr\'owka-Isbell spaces.
With the above remarks, it becomes clear that every Alexandrov-Urysohn compactum
 is separable, scattered and its third derivative set is empty.

Consider again the space $K=\ult(\fB)$ arising from some almost disjoint family as described above.
Recall that $M(\fB)$ is the dual of $C(\ult(\fB))$, see Section \ref{jls}.

\begin{lemma} \label{p:3}
Every  measure $\nu\in M(\fB)$ can be written as $\nu=\sum_{t\in K} c_t\delta_t$ where
$\sum_{t\in K}|c_t|=\|\nu\|<\infty$.
\end{lemma}

\begin{proof}
This follows simply from the fact  that every regular finite measure on a scattered compactum
is concentrated on a countable set,  see Rudin \cite[Theorem 6]{rudin57}.
\end{proof}

\begin{remark} \label{addrem}
If we write $\nu=\sum_{t\in K} c_t\delta_t$ as in Lemma \ref{p:3} then
the measure $\mu=\sum_{s\in S} c_s\delta_s$ will be called the $\ell_1(S)$-part of $\nu$.
Note that $\bar\nu=\nu-\mu$ is then a measure vanishing on all finite subsets of $S$.
\end{remark}

\section{The framework} \label{fr}
We now describe the basic features of the construction of the spaces leading to the proof of Theorem \ref{main}.
It will be convenient to work with the set  $\omega\times 2$ rather than with $\omega$.
Every subset of $\omega\times 2$ of the form $C=C_0\times 2$  will be called a \emph{cylinder}.
Moreover, given $n\in \omega$, we write $c_n$ for the cylinder $\{(n,0), (n,1)\}$.

Given a cylinder $C$, let us say that that a set $B\sub C$  {\em splits} $C$
if $|B\cap c_n|=1$ whenever $c_n\sub C$.
In such a case $C\sm B$ also splits $C$ and we say that the pair of sets $B$ and $C\sm B$ is a \emph{splitting}
of $C$.

We shall repeatedly consider the following framework:

\begin{enumerate}[FR(1)]
\item  $\cA$ is an almost disjoint family of cylinders in $\omega\times 2$;
\item for every cylinder $A\in \cA$, the pair $B_A^0, B_A^1$ is a splitting of $A$;
\item $\cB$ is an almost disjoint family $\{B_A^i\colon A\in\cA, i\in \{0,1\}\}$;
\item $\fA$ is the algebra of subsets of $\omega\times 2$ generated by
$\cA\cup \{c_n\colon n\in\omega\}$;
\item $\fB$ is the algebra generated by $\cB\cup \fin(\omega\times 2)$.
\end{enumerate}

Within this framework, we will slightly abuse notation concerning $JL$-spaces. Let us stress
that we declare $\jl(\cA)$ to be the closed subspace of $\ell_\infty(\omega\times 2)$
generated by characteristic functions of members in $\cA$, the constant one function and
characteristic functions of all finite \emph{cylinders} (rather than all finite sets).

The relatively simple result given below  is the first  basic ingredient of our proof of
  Theorem \ref{main}.

\begin{proposition}\label{stone2}\label{jl:2}
	In the setting FR(1) -- FR(5), the following hold:
	\begin{enumerate}
		\item[i)] There exists a commutative diagram
		\begin{center}
			\begin{tikzcd}
				\jl(\cA) \arrow[d, "i"] \arrow[r, "\sim"]  & C(\ult(\fA)) \arrow[d, "\theta^\circ"] \\
				\jl(\cB) \arrow[r, "\sim"] & C(\ult(\fB))
			\end{tikzcd}
		\end{center}
		\noindent where $i\colon\jl(\cA)\to\jl(\cB)$ is the inclusion map and
		$\theta^\circ\colon C(\ult(\fA))\to C(\ult(\fB))$
		is the isometric embedding $\theta^\circ (g) = g \circ \theta$ arising from the surjection $\theta$ as
in Lemma \ref{stone1}.
		
		\item[ii)] There is a closed subspace $X$ of  $\jl(\cB)$ such that $\jl(\cB)=\jl(\cA)\oplus X$ and
		both $X$ and $\jl(\cA)$ are  {\em 1}-complemented in $\jl(\cB)$.
	\end{enumerate}
\end{proposition}

\begin{proof}
	To prove $(i)$ it is enough to check that the diagram in question commutes
	for all $\chi_A$,  $A\in\fA$. Going down and then to the right we sent $\chi_A$ to
	the characteristic function of the clopen set
	\[C=\{p\in\ult(\fB)\colon A\in p\}\sub\ult(\fB),\]
	while
	the upper isometry  sends $\chi_A$ to the characteristic function of
	the clopen set
	\[D=\{q\in\ult(\fA)\colon A\in q\}\sub\ult(\fA).\]
	Clearly, $\chi_C=\chi_D\circ\theta$, so we are done.

\par Regarding $(ii)$, for every $g\in s(\fB)$ and $n\in\omega$ we define
\[ Pg(n,0)=Pg(n,1)=1/2\big(g(n,0)+g(n,1)\big).\]
Then $P$ is a norm one operator on $s(\fB)$ and $P\chi_{B_A^i}=1/2 \chi_A$ for every $A\in\cA$ and $i=0,1$
 so,
by Lemma \ref{jl:1}, $Pg\in \jl(\cA)$ for every $g\in s(\fB)$. It follows that
$P$ uniquely extends to a norm one projection from $\jl(\cB)$ onto $\jl(\cA)$.
We now consider $X = \ker P$ and a projection $Q = I-P$  onto $X$. We have
\[Qg(n,0) = -Qg(n,1) = 1/2\big(g(n,0)-g(n,1)\big), \]
so $\|Q\|=1$.
\end{proof}

\par

It is now clear that to prove  Theorem \ref{main}
we {\em only}  need to construct two almost disjoint families $\cA$ and $\cB$ satisfying FR(1) -- FR(5)
in such a way that the space $X$ defined in the proof of Proposition \ref{jl:2}(ii) is not a $C$-space.
This is, however, the hardest part of our  argument  requiring a number of preparations.

\subsection{More on measures from $M(\fB)$}
\label{fr:3}
Before embarking into our construction, we need some technical preparations concerning
$M(\fB)=\jl(\cB)^\ast$, where $\cB$ is an almost disjoint family of subsets of $\omega\times 2$ and $\fB = [\cB \cup \fin(\omega\times 2)]$, just as in
our framework FR(1) -- FR(5). We start with a simple remark.

\begin{remark}\label{jl:3}
	Consider the space  $X$ from  Proposition \ref{jl:2}(ii) and
	let $x^\ast\in X^\ast$. Then $x^\ast Q\in \jl(\cB)^\ast$. It follows that
	every functional on $X$ may be seen as a measure from $M_1(\fB)$ vanishing
	on all the cylinders from $\fB$.
\end{remark}

Our second observation provides a necessary condition for a subset of $M_1(\fB)$ to be free.

\begin{lemma}\label{p:2.5}
If a set of measures $M\sub M_1(\fB)$ is contained
in a free set, then for every $Z\in\fB$ and $\eps>0$ there
is  $g\in s(\fB)$ such that for every $\mu\in M$ we have
$\big| \bil{\mu}{g} -|\mu(Z)| \big| <\eps.$
\end{lemma}

\begin{proof}
As $|\mu(Z)|= |\bil{\mu}{\chi_Z}|$,
the function
$ \mu\longmapsto |\mu(Z)|$
 is weak$^\ast$ continuous. Therefore, since $M$
lies inside a free set, there is $f\in\jl(\cB)$ such that
$|\mu(Z)|=\bil{\mu}{f}$. To finish the proof it is enough to take
$g\in s(\fB)$ such that  $\|f-g\|<\eps$.
\end{proof}

Observe that we may identify $\ell_1(\omega\times 2)$ with a subspace of  $M(\cP(\omega\times 2))$:
every  element $x\in \ell_1(\omega\times 2)$ defines a measure $\mu$ by
 \[ \mu(B)=\sum_{s\in B} x(s), \quad B\sub\omega\times 2.\]
  In the sequel it will be customary to take some $\mu\in \ell_1(\omega\times 2) $
and consider the restriction of
$\mu$
to  $\fB$, expressing  such an intention by writing $\mu\in \ml(\fB)$.
By Remark \ref{addrem},
every measure $\nu\in M(\fB)$ can be written as $\nu=\mu+\bar\nu$, where
$\mu\in \ml(\fB)$ and $\bar\nu\in M(\fB)$ is a measure vanishing on finite sets and for which
$|\bar\nu|(B)\neq 0$ for at most countably many $B\in\cB$.

 Let us now consider $\cB_1\sub\cB$ and the corresponding subalgebra
$\fB_1=[\cB_1\cup\fin(\omega\times 2)] \subseteq \fB$. Every $\nu_1\in M(\fB_1)$
has an obvious extension to a measure $\nu\in M(\fB)$. Indeed, write
$\nu_1=\mu_1+\bar\nu_1$ where $\mu_1$ is the $\ell_1$-part of $\nu_1$ and extend $\bar\nu_1$
to a measure $\bar\nu\in M(\fB)$ by simply declaring that $\bar\nu(B)=0$ for every $B\in\cB\sm \cB_1$.
Then $\nu=\mu_1+\bar\nu\in M(\fB)$.

\begin{definition}\label{p:4}
In the above setting, the measure $\nu$ will be called the \emph{natural extension} of $\nu_1$.
\end{definition}

Finally, we include a suitable re-statement of the well-known Rosenthal lemma, see  e.g.\  \cite[p.18]{DU77}. Recall that $c_n=\{(n,0),(n,1)\}$.

\begin{theorem}\label{ros}
Let $(\mu_{n})_{n\in \omega}$ be a uniformly bounded sequence of measures in $\ell_1(\omega\times 2)$.
Then for every $\eps>0$ and every infinite set $N\sub\omega$ there is an infinite set $J\sub N$
 such that for every $k\in J$ we have
\[|\mu_k|\bigg( \bigcup_{j\in J\setminus\{k\}} c_j\bigg)<\eps.\]
\end{theorem}

\section{Separating sets of measures}\label{haydon}

We describe here a useful trick that can be used during an inductive construction of length $\con$; its essence
was invented by Haydon \cite[Proposition 3.5]{Ha81}.
Here the nature of the underlying countable set is inessential so we discuss measures on  $\omega$.
The main result of the section, Proposition \ref{haydon:prop} given below,
is stated  in quite a general form, for arbitrary measures on $\omega$.
In the next sections  it will be  applied to  $\ell_1$-measures,
 to prevent certain sets of measures from being free.

Throughout the section we fix a Boolean algebra $\fB\sub \cP(\omega)$ containing all finite subsets of $\omega$ and
such that $|\fB|<\con$.
We are concerned with pairs of sets of measures on $\omega$
which are ``separated with respect to $\fB$" in the sense of the following definition.

\begin{definition} \label{defi:sep}
	Consider two sets of measures $M,M'\sub M_1(\cP(\omega))$. We say that the pair $(M, M')$ is \emph{$\fB$-separated}
	if there
	are $\eps>0$, $n\in\omega$ and $B_1,\ldots, B_n\in\fB$ such that for every $(\mu, \mu')\in M\times M'$
	there is $i\le n$ such that $|\mu(B_i)-\mu'(B_i)|\ge \eps$.
\end{definition}

We start by recording  two useful properties of $\fB$-separation.

\begin{lemma}\label{haydon:1}
	Suppose that for given $M,M'\sub M_1(\cP(\omega))$ there are a simple $\fB$-measurable function
	$g\colon\omega\to\er$ and $\eps>0$ such that
	\[\big| \bil{\mu}{g}-\bil{\mu'}{g}\big| \ge\eps \mbox{ whenever } \mu\in M \mbox{ and }\mu'\in M'.\]
	Then the pair $(M,M')$ is $\fB$-separated.
\end{lemma}

\begin{proof}
	It follows by straightforward calculations. Indeed, if we write
	$g=\sum_{i=1}^n r_i\cdot\chi_{B_i}$ for some $r_i\in\er$ and $B_i\in\fB$,
	then the pair $(M,M')$ is $\fB$-separated by the sets $B_i$, $1\leq i\leq n$,  with $\eps' = \eps/ \left(\sum_{i=1}^n |r_i|\right)$ as the separating constant.
\end{proof}

\begin{lemma}\label{haydon:2}
	For every infinite subset $M\sub M_1(\cP(\omega))$
	the pair $(M',M\sm M')$ can be $\fB$-separated
	for fewer than $\con$ many sets $M'\sub M$.
\end{lemma}

\begin{proof}
	Suppose that $M'\sub M \sub M_1(\cP(\omega))$ and the pair $(M',M\sm M')$ is separated by
	some $B_1, ..., B_n\in \fB$ and some $\eps>0$. Choose a rational number $\delta<\eps/2$ and for any
	$q:=(q^1, ..., q^n)\in [-1,1]^n\cap {\mathbb Q}^n$ consider the weak$^\ast$ open set
	\[ M(q)=\{\mu\in M_1(\cP(\omega))\colon |\mu(B_i)-q^i|<\delta\mbox{ for every } i\le n\}.\]
	Note that, by separation, if $M(q)\cap M'\neq\emptyset$ then $M(q)\cap (M\sm M')=\emptyset$.
	The  compactness of $[-1,1]^n$ allows us to obtain a finite cover of $M_1(\cP(\omega))$ consisting of sets
	of the form $M(q)$ for $q\in [-1,1]^n \cap \mathbb Q^n$. Therefore, we deduce the existence of $m\in \omega$ and $q_1, ..., q_m \in [-1,1]^n \cap {\mathbb Q}^n$ such  that
	\[ M' \subseteq \bigcup_{1\leq k \leq m} M(q_k)\quad\mbox{and}\quad
	(M\sm M') \cap \bigcup_{1\leq k \leq m} M(q_k) = \emptyset.\]
	In other words, $M'$ and $M\sm M'$ can be ``physically separated" by some finite union of sets
	of the form $M(q)$. Finally, it is clear that there are less than $\con$ such unions, and each of them
	can separate only one pair $(M',M\sm M')$.
\end{proof}

Given any $Z\sub\omega$, we shall write $\fB[Z]$ for the algebra generated by $\fB\cup\{Z\}$.
Note that every element of $\fB[Z]$ is of the form $(A\cap Z)\cup(B\cap Z^c)$ for some $A,B\in\fB$.
We are finally ready for the main result of the section, the following
measure-theoretic version of Haydon's Lemma \cite[1D]{Ha81}.

\begin{proposition}\label{haydon:prop}
	Let us fix $\xi<\con$. Suppose we are given
	\begin{enumerate}
		\item[i)] a list $\{\fB_\alpha\colon \alpha<\xi\}$ of Boolean subalgebras of $\cP(\omega)$ such that,
		for every $\alpha<\con$, $\fin(\omega) \subseteq \fB_\alpha$ and $|\fB_\alpha|\le |\xi|$;
		
		\item[ii)] a list $\{(M_\alpha, M'_\alpha)\colon \alpha<\xi\}$ of pairs of subsets of $M_1(\cP(\omega))$ such that, for every $\alpha<\con$, the pair $(M_\alpha, M'_\alpha)$ is not $\fB_\alpha$-separated.
	\end{enumerate}
	Then for every almost disjoint family $\cZ\sub \cP(\omega)$ with $|\cZ|=\con$
	there exists $Z\in\cZ$ such that for every $\alpha<\xi$, the pair $(M_\alpha, M_\alpha')$ is not
	$\fB_\alpha[Z]$-separated.
\end{proposition}

The proof of Proposition \ref{haydon:prop} builds on the following technical lemma.

\begin{lemma}\label{haydon:3}
	Suppose that for given $M,M'\sub M_1(\cP(\omega))$ there are
	
	\begin{itemize}
		\item $n\in\omega$, $\eps>0$ and $k>12n/\eps$;
		\item almost disjoint sets $Z_1,\ldots, Z_k\sub\omega$;
		\item $A_1,B_1, ..., A_n, B_n\in \fB$;
	\end{itemize}
	such that, for every $1\le j \leq k$, the pair $(M,M')$ is $\fB[Z_j]$-separated, with $\varepsilon$ as a constant, by
	the sets
	\[ Y_{j,i}=(A_i\cap Z_j)\cup (B_i\cap Z_j^c),\quad 1\leq i \leq n.\]
	Then $(M,M')$ is $\fB$-separated.
\end{lemma}

\begin{proof} Let us first prove the lemma. Consider a fixed pair $(\mu,\mu')\in M\times M'$. For every $j\le k$,  find $i(j)\le n$
	such that $|\mu(Y_{j,i(j)}) -\mu'(Y_{j,i(j)})|\ge\eps$. Since $k>12n/\eps$,
	we can choose $m\le n$ such that the set of indices $J=\{j\le k\colon i(j)=m\} $ verifies $|J|> 12/\eps$.
	Now set \[F=\bigcup\{Z_j\cap Z_l\colon j,l\in J, j\neq l\}.\]
	It is clear that, since $Z_1,\ldots, Z_k$ are almost disjoint, $F$ is finite and the sets $Z_j\sm F$, where $j\in J$, are pairwise disjoint. Consequently, there is a privileged  $j\in J$ such that
	\begin{equation} \label{eq:j} |\mu|(Z_j\sm F)<\eps/6 \quad \mbox{and} \quad |\mu'|(Z_j\sm F)<\eps/6.
	\end{equation}
	Indeed, since $\|\mu\|\le 1$, writing $J_1$ for the set of those
	$j\in J$ for which    $|\mu|(Z_j\sm F)\ge\eps/6$ we have
	\[1\ge \|\mu\|\ge \sum_{j\in J_1}|\mu|(Z_j\sm F)\ge |J_1|\cdot\eps/6,\]
	so $|J_1|\le 6/\eps<|J|/2$. Of course, the same argument applies to $\mu'$, and this
	shows that some $j\in J$ satisfies (\ref{eq:j}).
	
	We take $m$ as above and $j$ satisfying (\ref{eq:j}). Put
	\[ Y= Y_{j,m}=(A_m\cap Z_j)\cup (B_m\cap Z_j^c)\]
	and note that we have  $|\mu(Y)-\mu'(Y)|\ge \eps$.
	We now define
	\begin{equation} \label{eq:b} C = (A_m\cap Z_j\cap F)\cup (B_m\cap (Z_j\cap F)^c).
	\end{equation}
	Observe that, since $F$ is finite, $C\in \fB$. Also, it is not difficult to see that
	$C\sm Y$ and $Y\sm C$ are both contained in $Z_j\sm F$, just by observing that
	for $x\notin Z_j\sm F$ we have $x\in Z_j$ if and only if $x\in Z_j\cap F$.
	%
	Straightforward calculations using (\ref{eq:j}) yield
	\[ |\mu(Y)-\mu(C)|< \eps/3 \quad \mbox{ and } \quad |\mu'(Y)-\mu'(C)|< \eps/3,\]
	and so we finally deduce that $ |\mu(C)-\mu'(C)|\ge \eps/3$.
	\par To conclude the proof of Lemma \ref{haydon:3}, it is enough to realise that there are finitely many sets $C$ defined as in (\ref{eq:b}),
	and therefore the pair $(M,M')$ is $\fB$-separated by the collection of all such sets with $\eps/3$ as the separating constant.
	\medskip

	\par We can now prove
Proposition \ref{haydon:prop} by a simple counting argument, following \cite[Proposition 3.5]{Ha81}.
Suppose that every $Z\in\cZ$ fails to satisfy the thesis of Proposition \ref{haydon:prop}, and this is witnessed by some
$\alpha<\xi$, some $n\in\omega$, some rational $\eps>0$
and some collection $A_i, B_i\in\fB_\alpha$ for $i\le n$. Since the number of such witnesses is less than $\con$, there are infinitely many
$Z\in\cZ$ having the same witness; then Lemma \ref{haydon:3} implies that
the pair $(M_\alpha, M_\alpha')$ is $\fB_\alpha$-separated, contrary to our assumption.
\end{proof}

\section{Separation versus freeness}
\label{sec:alg}

In this section we fix an infinite cylinder  $C=C_0\times 2 \sub\omega\times 2$.
Our objective is to describe some method of  defining a {\em new} cylinder $C'\sub  C$ and conveniently split it into two infinite subsets $B^0$ and $B^1$. We shall  say that a Boolean algebra $\fB \subseteq \cP(\omega \times 2)$ which contains $\fin(\omega\times 2)$ is \emph{trivial on} $C$ if either $B\cap C$ or $B^c\cap C$ is finite
for every $B\in\fB$.

\par Let us recall that, in the context of our framework FR(1) -- FR(5), the space $X$ in the direct sum
$\jl(\cB) = \jl(\cA)\oplus X$ contains every function of the form
$f_n = \frac12(e_{(n,0)} - e_{(n,1)})$ for $n\in \omega$.
Therefore, every norming subset for $X$, when regarded as a subset of $M(\fB)$
as indicated in Lemma \ref{p:3},  must contain a sequence  $(\nu_n)_{n\in\omega}$ such that $\inf_n |\nu_n(f_n)| = \inf_n |\nu_n(n,0)|>0$. This, together with the description of elements of $X^*$ as measures in $M(\fB)$ vanishing on all the cylinders of $\fB$ --- see Remark \ref{jl:3} --- motivates the following working definition.

\begin{definition} \label{admissible} We say that a sequence
	$(\mu_n)_{n\in \omega}$ in the unit ball of $\ell_1(\omega\times 2)$  is \emph{admissible} if
	$\mu_n(c_k)=0$ for all $k,n\in \omega$ and		 $\inf_n |\mu_n(\{n,0\})|>0$.
\end{definition}

It will be also convenient to use the following piece of notation:
given a sequence $(\mu_n)_{n\in\omega}$ and any set $J\sub\omega$, we write
\[\mu[J]=\{\mu_n: n\in J\}.\]

\begin{remark}\label{additional_remark}
	Suppose that $\fB$ is trivial on the cylinder $C$ and consider an admissible sequence
	$(\mu_n)_{n\in \omega}$. Then for any $J\sub I\sub\omega$, the pair
	$\big(\mu[J],\mu[I\sm J]\big)$ is $\fB[C]$-separated if and only if it is $\fB$-separated.
	Indeed, the key point is that for every $D\in \fB[C]$
	there is $B\in\fB$ such that $\mu_n(D)=\mu_n(B)$ for every $n\in\omega$. To prove this, observe that  $\mu_n(C)=\mu_n(\omega\times 2)=0$, by the first requirement of admissibility. Now, if $D=(B_1\cap C)\cup (B_2\cap C^c)$ for suitable $B_1, B_2\in \fB$, and $B_1\cap C$ is not finite, then necessarily $B_1^c \cap C \in \fB$ and $\mu_n(B_1\cap C) = \mu_n((B_1^c\cap C)^c)$. A similar argument can be applied to $B_2\cap C^c$.
\end{remark}

\begin{lemma} \label{lemma:separation} Let $(\mu_n)_{n\in\omega}$
	be an admissible sequence.
	There are three infinite sets $ J_2\sub J_1\sub J_0\subseteq C_0$ and a set
	$Z\subseteq \omega \times 2$ that splits the cylinder $ J_1\times 2 $
	such that for every algebra $\fB$ that  is trivial on $C$,
	if the sequence $(\mu_n)_{n\in\omega}$ is inside
	a free subset of $M_1(\fB[Z])$, then at least of one of the pairs
	\[ \big(\mu[J_2], \mu[ J_1\setminus J_2]\big), \ \big(\mu[ J_1] , \mu[J_0\setminus J_1]\big),\]
	is $\fB$-separated.
\end{lemma}

\begin{proof} First, let $c=\inf_n |\mu_n(\{n,0\})|$, and pick $0<\delta<c/16$. By virtue of Rosenthal's Lemma (see Theorem \ref{ros}), there is an infinite $J_0\sub C_0$
	such  that
	\begin{equation} \label{eq:1}
		|\mu_n|\big(( J_0\times 2)  \setminus c_n)\big) <\delta \mbox{ for every }n\in J_0.
	\end{equation}
	
	Then  find an infinite subset $J_1\subseteq J_0$ with $J_0\sm J_1$ infinite,
	and singletons $p_n \sub c_n$ so that the sequence $(\mu_n(p_n))_{n\in J_1}$ is convergent to some $a\ge c$.
	By eliminating a finite number of elements from $J_1$, we can assume that
	\begin{equation} \label{eq:2}
		|\mu_n(p_n) - a|<\delta \mbox{ for every } n\in J_1.
	\end{equation}
	Pick now any infinite set $J_2\sub J_1$ such that $J_1\sm J_2$ is also infinite and put	
	\[ Z = \big(\bigcup_{n\in J_2} p_n \big) \cup \big(\bigcup_{n\in J_1 \setminus J_2} (c_n\setminus p_n)\big).
	\]
	
	To continue, we first introduce a bit of notation:
	given $r,s\in \R$ and $\eps>0$, it will be convenient to write $r\approx_\eps s$ rather than $|r-s|<\eps$.
	
	Consider an algebra $\fB$ as in the assertion of Lemma \ref{lemma:separation} and
	let us suppose
	that the sequence $(\mu_n)_{n\in \omega}$ lies inside a free set  $M\sub M_1(\fB [Z])$.
	Hence Lemma \ref{p:2.5} ensures that
	there is a $\fB [Z]$-measurable simple function $h\colon \omega\times 2 \to \R$ such that
	\[ |\mu_n(Z)|\approx_\delta \bil{\mu_n}{h} \text{ for } n\in \omega.\]
	
	For our purposes there is no loss of generality if we assume that $h=r\cdot \chi_Z+g$, where $r\in \R$ and
	$g$ is $\fB$-measurable. Indeed, as $\fB$ is trivial on $J_0\times 2$ and contains all finite sets,
	we can additionally assume that $J_0\times 2\in \fB$, see Remark \ref{additional_remark}.
	Consequently,  every $\fB[Z]$-measurable simple
	function can be written as
	\[r_1\cdot \chi_Z+r_2\cdot\chi_{(J_0\times 2)\sm Z}+g=(r_1-r_2)\cdot \chi_Z + (r_2\cdot\chi_{J_0\times 2}+g),\]
	where $g$ is $\fB$-measurable and so is the second summand (as $J_0\times 2\in\fB$).
	
	Let us further suppose that $r\geq 0$
	(otherwise, we may consider the  function $\mu \longmapsto -|\mu(Z)|$).
	Then we have
	\[ |\mu_n(Z)| \approx_\delta r\cdot\mu_n(Z) + \bil{\mu_n}{g} \text{ for every } n\in \omega. \]
	Using equations (\ref{eq:1}) and (\ref{eq:2})  we conclude that
	\begin{equation} \label{eq:3} \left\{ \begin{array}{lcrl}
			a & \approx_{\delta(3+2r)} &  r\cdot a + \bil{\mu_n}{g} & \mbox{ for  }n\in J_2, \\[2mm]
			a & \approx_{\delta(3+2r)} & -r\cdot a +  \bil{\mu_k}{g} & \mbox{ for } k\in J_1 \setminus J_2,  \\[2mm]
			0 & \approx_{\delta(3+2r)}  &  \bil{\mu_l}{g}& \mbox{ for } l\in J_0 \setminus J_1.
		\end{array}  \right. \end{equation}
	
	Now we plan to deduce that the above relations (\ref{eq:3}) imply separation of one of the pairs from
	the conclusion of Lemma \ref{lemma:separation}.
	Suppose that $r$ is small, say $0\le r<1/2$. Then the first two relations in (\ref{eq:3}) give
	\[\bil{\mu_n}{g}\ge \eps_1 :=a/2-\delta(3+2r)\ge a/2- 4\delta \mbox{ for all }n\in J_1,\]
	while the third one 	gives
	\[ \bil{\mu_l}{g}\le \eps_2 :=\delta(3+2r)\le 4\delta \mbox{  for  }l\in J_0 \setminus J_1.\]
	Since $\eps_1>\eps_2$	(recall that $a\ge c$ and $\delta<c/16$),
	Lemma \ref{haydon:1} implies  that
	$\mu[ J_1]$ and $\mu[ J_0 \setminus J_1]$
	are $\fB$-separated.
	
	On the other hand, if $r\ge 1/2$, then for $n\in J_2$ and $k\in J_1\sm J_2$ the first two
	approximate relations in (\ref{eq:3}) yield
	\[ \bil{\mu_k}{g}-  \bil{\mu_n}{g} \ge  2ra-2\delta(3+2r)
	=2r(a-2\delta)-6\delta \ge a-8\delta>0,\]
	so this time
	$\mu[ J_2]$ and $\mu[ J_1\setminus J_2]$ are $\fB$-separated, again by Lemma \ref{haydon:1}.
\end{proof}

\begin{remark}\label{lemma:separation2}
	Observe that the set $J_1$ defined in the proof above can be replaced by any other infinite subset of $J_1$. The same condition is true for the set $J_2$, which is only required to satisfy that both $J_2$ and $J_1\sm J_2$ are infinite.
\end{remark}

For the purpose of the final construction it will be convenient to augment slightly
the previous lemma as follows.

\begin{lemma} \label{lemma:separation3}
	In the setting of Lemma \ref{lemma:separation},  the sets $J_0, J_1,J_2$ and $Z$ constructed in its proof satisfy the following assertion.
	
	Consider an algebra $\fB$ that   is trivial on $C$, and
	a sequence $(\bar\nu_n)_{n\in\omega}$ of measures on  $\fB[J_0\times 2]$ which vanish on finite sets and are identically zero on $J_0\times 2$. Then, writing  $\nu_n=\mu_n+\bar\nu_n$,
	if the sequence of  $\nu_n$  is inside a
	free subset of $M_1(\fB[Z])$, then one of  the pairs
	$\big(\nu[ J_2], \nu[J_1\setminus J_2]\big)$,
	$\big(\nu[J_1], \nu[J_0\setminus J_1]\big)$
	is $\fB$-separated.
\end{lemma}

\begin{proof}
	This follows by the same argument as above:
	as every $\bar\nu_n$ is identically zero
	on $J_0\times 2$,
	all the considerations inside $J_0\times 2$ do not affect the sequence of measures $(\bar\nu_n)_{n\in \omega}$.
\end{proof}

\section{The final construction}
\label{sec:final}
We are ready to present the proof of our main result, Theorem \ref{main}, by constructing
an almost disjoint family $\cA=\{A_\alpha\colon\alpha<\con\}$ of cylinders in $\omega\times 2$  and
defining a splitting of every $A_\alpha$ into two sets $B_\alpha^0, B_\alpha^1$ in accordance with our
framework FR(1) -- FR(5). The main point is to prevent
the space $X$ of Proposition \ref{jl:2} from being  a $C$-space. The construction
will be performed by induction  of length $\con$.

\par Let us first introduce the following notation. For every $V\sub\con$ we write
\[
\fB(V)=\big[\{B_\alpha^0, B_\alpha^1\colon \alpha\in V\}\cup\fin(\omega\times 2)\big].\]
Note that $\fB(\alpha)$ stands for $\fB(\{\beta\colon \beta<\alpha\})$, and
in particular, $\fB(0)=\big[\fin(\omega\times 2)\big]$.
Our final algebra will be $\fB:=\fB(\con)$.

It is convenient to realize that, even before we start our construction, we can code all the measures from $M_1(\fB)$
in the following way.
Every $\nu\in M_1(\fB)$ can be written as  $\nu=\mu+\bar\nu$ where $\mu$ is its
$\ell_1(\omega\times 2)$-part as in Remark \ref{addrem}. Then $\mu$ is uniquely determined
by its values on singletons and $\bar\nu$ is uniquely determined by the values of
$\bar\nu(B_\alpha^i)$ for $\alpha<\con$, $i=0,1$.
Recall that, in view of Remark \ref{jl:3}, we are only interested in those
$\nu\in M_1(\fB)$ that vanish on all the cylinders from $\fB$.
Let $\Psi$ be the set of all $(x,y)\in \ell_1(\omega\times 2)\times \ell_1(\con\times 2)$
such that

\begin{enumerate}[--]
\item $x(n,0)+x(n,1)=0$ for every $n\in\omega$ and $y(\alpha,0)+y(\alpha,1)=0$ for every $\alpha<\con$;
\item $\|x\|_1+\|y\|_1\le 1$.
\end{enumerate}

Then the set $\Psi^\omega$ codes all the sequences of such measures;
however, during our construction we only need to treat codes of sequences of measures
that appear in norming subsets of $X^*$. Let $\Phi\sub\Psi^\omega$ be defined by the condition
\[\vf=(\psi_n)_{n\in\omega}\in\Phi \mbox{ if, writing }  \psi_n=(x_n, y_n) \mbox{, we have } \inf_n |x_n(n,0)|>0.\]
Let us enumerate $\Phi$ as $\{\vf_\alpha\colon \alpha<\con\}$ in such a way that
\[\mbox{ if } \vf_\alpha=((x_n)_n,(y_n)_n)\in  \Phi \mbox{ then }
y_n(\xi,i)=0 \mbox{  for  }i=0,1  \mbox{ and every } \xi\ge \alpha.\]

The enumeration mentioned above can be obtained using the fact that
for every $y\in\ell_1(\con\times 2)$ there is $\alpha<\con$ such that
$y$ is zero outside $\alpha\times 2$. For instance, take first
an enumeration $\{\vf_\alpha'\colon \alpha<\con\}$ of $\Phi$ in which every of its elements
is present cofinally often. Then for $\vf_\alpha'=(x_n,y'_n)_n$ put
$\vf_\alpha=(x_n,y_n)_n$ where, for every $n\in\omega$, $y_n(\beta)=y_n'(\beta)$ whenever $\beta<\alpha$
and $y_n(\beta)=0$ for $\beta\ge\alpha$.

At the very beginning we  fix an almost disjoint family $\cR = \{R_\alpha: \alpha<\con\}$ of infinite subsets of $\omega$.
The role of $R_\alpha$ is to make room for step $\alpha$.

In the sequel, we write $(\nu_n^\alpha)_{n}$ to refer to the sequence in
 $M(\fB)$ coded by $\varphi_\alpha$.
For every $\alpha<\con$ we construct three infinite subsets
$ J^\alpha_2\sub J^\alpha_1\sub J^\alpha_0 \sub R_\alpha$ and  a splitting
 $B_\alpha^0, B_\alpha^1$ of the cylinder
$J^\alpha_1 \times 2$ so that the following is satisfied for all $\xi<\con$:

\begin{ind} For every $\alpha < \xi$, the pairs
	\[ \big(\nu^\alpha[ J^\alpha_2], \nu^\alpha[ J^\alpha_1\sm J^\alpha_2]\big)
	\mbox{ and  } \big(\nu^\alpha[ J^\alpha_1], \nu^\alpha[ J^\alpha_0\sm J^\alpha_1]\big) \]
	are not separated by the algebra $\fB(\xi\sm\{\alpha\})$.
\end{ind}

\newcommand{\kp}{\protect{\bf KP}}

Note that $\kp(0)$ is  vacuously satisfied and if $\kp(\eta)$ holds for every $\eta$ below some limit
ordinal $\xi$ then $\kp(\xi)$ is granted (an increasing family of `non-separating' algebras
is again `non-separating').

 Assuming that \kp($\xi$) holds, we need to show how to perform the successor step and guarantee \kp$(\xi+1)$.
Given the sequence $(\nu^\xi_n)_{n\in\omega}$ coded by  $\vf_\xi$, we write
$\nu^\xi_n = \mu^\xi_n + \bar\nu^\xi_n$, where $\mu^\xi_n$ is the $\ell_1(\omega \times 2)$-part
(see Remark \ref{addrem}).
We now fix an almost disjoint family $\cC$ of size $\con$ consisting of infinite cylinders contained in
$R_\xi\times 2$.
We apply Lemma \ref{lemma:separation3}  to the algebra $\fB(\xi)$ $\con$-many times,
each time for a given  cylinder $\cC\ni C=C_0\times 2$. We  get the  sets $J_2(C), J_1(C), J_0(C)\sub C_0$ and
$Z(C)\sub J_1(C)\times 2$ satisfying
the assertion of \ref{lemma:separation3}.
Using  Lemma \ref{haydon:2} in combination with Remark \ref{lemma:separation2},
we can assume that both pairs
$ \big(\nu^\xi[ J_2(C)], \nu^\xi[ J_1(C)\sm J_2(C)]\big)$ and
$ \big(\nu^\xi[J_1(C)], \nu^\xi[ J_0(C) \sm J_1(C)]\big)$ are not $\fB(\xi)$-separated.

 Now Proposition \ref{haydon:prop} comes into play: we apply it
the family of algebras $\fB_\alpha=\fB(\xi\sm\{\alpha\})$, where $\alpha<\xi$, and
to the almost disjoint family  $\cZ = \{Z(C)\colon C\in\cC\}$.
The proposition enables us to pick
a suitable $C\in\cC$  having the property that  for every $\alpha<\xi$, the pairs
\[ \big(\nu^\alpha[ J^\alpha_2], \nu^\alpha[ J^\alpha_1\sm J^\alpha_2]\big)
\mbox{ and  } \big(\nu^\alpha[ J^\alpha_1], \nu^\alpha[ J^\alpha_0\sm J^\alpha_1]\big)\]
are  not separated even by the enlarged algebra $\fB(\xi\sm\{\alpha\})[Z(C)].$
Then we set
\begin{itemize}
	\item $A_\xi=J_1(C)\times 2,  B_\xi^0=Z(C), B_\xi^1= A_\xi \sm B_\xi^0$;
	\item $J_2^\xi=J_2(C),
	J_1^\xi=J_1(C), J_0^\xi=J_0(C)$.
\end{itemize}
By our choice,  \kp($\xi+1$) holds and this concludes the successor step.

Once the construction has been done, we are ready to complete the proof of Theorem \ref{main}.
The almost disjoint families $\cA=\{A_\alpha\colon \alpha<\con\}$ and
$\cB= \{B_\alpha^0, B_\alpha^1\colon \alpha<\con\}$
clearly satisfy the conditions FR(1) -- FR(5) of the framework. Hence
Proposition \ref{stone2} explains everything except for  the fact
that the space $X$ present in the direct sum $\jl(\cB) = \jl(\cA) \oplus X$ is not a $C$-space, which we need to
demonstrate.

Consider a norming set in $B_{X^\ast}$, and regard it as a subset of measures in
$M_1(\fB)$ vanishing on all the cylinders from $\fB$, see Remark \ref{jl:3}.
Then our norming set must contain a sequence of measures $(\nu_n)_{n\in\omega}$ such that
their $\ell_1(\omega \times 2)$-parts form an admissible sequence.
Consequently, $(\nu_n)_{n\in\omega}=(\nu_n^\xi)_{n\in\omega}$,
for some $\xi<\con$, i.e.\ the sequence in question
was coded by $\vf_\xi$.

We took care of the sequence $(\nu_n^\xi)_{n}$
at step $\xi$ when  we added elements $B_\xi^i$ to our algebra.
Since  \kp$(\eta)$ holds for every $\eta<\con$,
 $ \big(\nu[ J^\xi_2], \nu[ J^\xi_1\sm J^\xi_2]\big)
\mbox{ and  } \big(\nu[ J^\xi_1], \nu[ J^\xi_0\sm J^\xi_1]\big) $
are not $\fB(\con\sm\{\xi\})$-separated.
Hence Lemma \ref{lemma:separation3} says that
$(\nu_n^\xi)_{n}$ cannot lie  inaide  a free subset of $M_1(\fB)$.

Finally, it follows from Proposition \ref{p:1} that $X$ is not a $C$-space and
the proof is complete.

\section{Applications to twisted sums of $c_0$ and $C$-spaces}\label{section:twisted}

 Let us recall that an \emph{exact sequence} of Banach spaces is a diagram of Banach spaces and operators
\begin{equation*}  \begin{tikzcd}
	0 \arrow[r] & Y \arrow[r, "i"] & Z \arrow[r, "q"] & X \arrow[r] & 0
	\end{tikzcd} \end{equation*}
in which the kernel of every arrow coincides with the image of the preceding one. The middle space $Z$ is usually called a \emph{twisted sum} of $Y$ and $X$. By the open mapping theorem, $i(Y)$ is a closed subspace of $Z$ so that $Z/i(Y)$ is isomorphic to $X$. We also say an exact sequence is \emph{trivial}, or that \emph{it splits}, when $i(Y)$ is a complemented subspace of $Z$. The reader can find a fully detailed exposition on the topic in \cite{hmbst} and
 a discussion of twisted sums of $c_0$ and $C$-spaces in \cite{AMP20}.

There are two natural ways of constructing twisted sums of $C$-spaces: one is to consider
a surjection $\theta\colon L \to K$ between compact spaces
which produces the exact sequence
\begin{equation} \label{eq:theta-seq} \begin{tikzcd}
	0 \arrow[r] & C(K) \arrow[r, "\theta^\circ"] & C(L) \arrow[r] & X \arrow[r] & 0.
\end{tikzcd} \end{equation}
The other is the ``dual" situation: we consider an embedding $\iota\colon K \to L$, which gives
\begin{equation} \label{eq:iota-seq} \begin{tikzcd}
	0 \arrow[r] & \ker \iota^\circ \arrow[r] & C(L) \arrow[r, "\iota^\circ"] & C(K) \arrow[r] & 0.
\end{tikzcd} \end{equation}

\subsection{Deeper into Pe\l czy\'nski's problem}
\par It is not difficult to show that, in the setting of equation (\ref{eq:iota-seq}), $\ker \iota^\circ$ is always a $C$-space, since it is isomorphic to the space of continuous functions on the quotient space $L/\iota(K)$ obtained by identifying the subspace $\iota(K) \subseteq L$ to a point. The corresponding situation for sequences arising from continuous surjections, namely, if the quotient space of (\ref{eq:theta-seq}) is always a $C$-space, was the content of Problem \ref{quotient}, and Section \ref{sec:final} contains a negative answer. We now show a somewhat simpler example illustrating this phenomenon:

\begin{theorem} \label{additional}
There is a compact space $L$ and a short exact sequence
\[ \begin{tikzcd}
	0 \arrow[r] & c_0 \arrow[r] & C(L) \arrow[r] & W \arrow[r] & 0,
\end{tikzcd} \]
in which $W$ is not a $C$-space.
\end{theorem}

\begin{proof}
Let $\fB$ be the subalgebra of $P(\omega \times 2)$ constructed in Section \ref{sec:final} and write again
 $L=\ult(\fB)$. Consider
the subalgebra $\fA_0$ of $\fB$  generated by all the cylinders $c_n$.
Of course, $\fA_0$ is isomorphic to the Boolean algebra generated by all finite sets in $\omega$;
therefore $C(\ult(\fA_0))$ is the classical Banach space $c$ (of all converging sequences), which is
isomorphic to its hyperplane $c_0$.
If we write  $p\colon \ult(\fB) \to \ult(\fA_0)$ for the canonical quotient map then a very similar argument to
the one from the proof of our main theorem shows  that the quotient space $W$ in the sequence
\[ \begin{tikzcd}
		0 \arrow[r] & C(\ult(\fA_0)) \arrow[r, "p^\circ"] & C(L) \arrow[r] & W \arrow[r] & 0,
	\end{tikzcd} \]
is not a $C$-space. Indeed, the crucial point is that it is still true that every functional over $W$ can be identified with a finitely additive measure on $\fB$ vanishing on all the cylinders $C_n$. Therefore,
if $F$ is a $c$-norming set in $B_{X^\ast}$ then it must contain a sequence of measures in $M(\fB)$ so that their respective $\ell_1$-parts form a $c$-admissible sequence.  We now proceed exactly as in the proof of
Theorem \ref{main} to deduce that $F$ cannot be free, and so $W$ cannot be  a $C$-space.
\end{proof}

Note however, that the quotient space $W$ in Theorem \ref{additional}
 is not isomorphic to a complemented subspace of $C(L)$.
  This follows from the fact that $W$ contains a complemented copy of $c_0$; say,
  $W\simeq c_0 \oplus Y$.
  Then $C(L) \simeq W \oplus c_0$ would lead to
 \[W\simeq c_0\oplus Y\simeq c_0\oplus (c_0\oplus Y)\simeq c_0 \oplus W\simeq C(L),\]
so $W$ would be a $C$-space.

As for the mentioned copy of $c_0$ inside $W$, let us consider any $\xi < \con$ and recall that the cylinder $A_\xi = J^\xi_2 \times 2$ is split into $B_\xi^0$ and $B_\xi^1$. This allows us to write $B_\xi^0 = \{(n, k_0(n))\colon n\in J^\xi_2 \}$ and $B_\xi^1 = \{(n, k_1(n))\colon n\in J_2^\xi\}$ where $k_0(n) \neq k_1(n)$ for every $n\in J_2^\xi$. Finally, the subspace of $W$ spanned by the functions
\[ g_n = \chi_{(n, k_0(n))} - \chi_{(n, k_1(n))} , \quad n\in J_2^\xi \]
is isometrically isomorphic to $c_0$ and complemented in $W$. The projection is simply $P(f) = \sum_{n\in\omega} (f(n, k_0(n)) - f(p_\xi^0)\big) \cdot g_n$, where $p_\xi^0$ is the only point in $\ult(\fB)$ containing $B_\xi^0$ and no finite set, see Section \ref{sec:au}.

\subsection{Twisted sums of $c_0$ and $c_0(\con)$ which are not $C$-spaces}
\par In \cite{PS21} we showed, using some set-theoretic assumption,  that for
every Eberlein compact space $K$ of weight $\con$ there is a short exact sequence
\[ \begin{tikzcd}
	0 \arrow[r] & c_0 \arrow[r] & Z \arrow[r] & C(K) \arrow[r] & 0,
\end{tikzcd} \]
where $Z$ is not  a $C$-space. This construction was improved in \cite[Theorem 3.4]{CS22}, where the authors show that it remains true for other well-known classes of compact spaces. We now show that a twisted sum of $c_0$ and $c_0(\con)$ which is not a $C$-space can be obtained in ZFC. Indeed,  the space $X = \jl(\cB)/\theta^\circ[\jl(\cA)]$ constructed in Section  \ref{sec:final} is also a twisted sum of $c_0$ and $c_0(\con)$.

\par To show this, let us write $Q\colon C(L) \to X$ for the quotient map obtained in Proposition \ref{jl:2}(ii). We shall denote by $X_0$ the canonical copy of $c_0$ inside $\jl(\cB)$.
 It is not difficult to show that $Q(X_0)$ is isomorphic to $c_0$, since it is spanned by the sequence of functions
\[ f_{n} = Q(\chi_{(n,0)}) = \tfrac12(\chi_{(n,0)} - \chi_{(n,1)}), \quad n\in \omega.\]
To describe the quotient $X/Q(X_0)$, we consider the functions
\[ f_\xi = Q(\chi_{B_\xi^0}) = \tfrac12\big(\chi_{B_\xi^0} - \chi_{B_\xi^1}\big), \quad \xi < \con \]
and write $R\colon X \to X/Q(X_0)$ for the canonical quotient map. Then it is straightforward to see that $ X/Q(X_0)$ is spanned by $\{R(f_\xi)\colon \xi < \con\}$ and that for every $n\in \omega$, every $\lambda_i\in \R$ and $\xi_i<\con$ , $i\leq n$ we have $\|\sum_{i\leq n} \lambda_i R(f_{\xi_i})\| = \frac12\max_{i\leq n} |\lambda_i|$, thanks to the fact that $\cB$ is almost disjoint. This shows that $X/Q(X_0)$ is isomorphic to $c_0(\con)$.

\end{document}